\documentclass[12pt]{amsart}
\usepackage{amsfonts,amssymb, amscd, latexsym, graphicx, psfrag, color,float}

\usepackage[all]{xy}

\usepackage{booktabs}
\usepackage{multirow}
\usepackage{mathrsfs}

\usepackage[margin=1in,marginparwidth=0.8in, marginparsep=0.1in]{geometry}

\usepackage[bookmarks=true, bookmarksopen=true,%
bookmarksdepth=3,bookmarksopenlevel=2,%
colorlinks=true,%
linkcolor=blue,%
citecolor=blue,%
filecolor=blue,%
menucolor=blue,%
urlcolor=blue]{hyperref}

\swapnumbers

\theoremstyle{definition}
\newtheorem*{prop*}{Proposition}

\newcommand{\bF}{\mathbf{F}}
\newcommand{\bG}{\mathbf{G}}
\newcommand{\bQ}{\mathbf{Q}}
\newcommand{\bR}{\mathbf{R}}
\newcommand{\bZ}{\mathbf{Z}}

\newcommand{\Hom}{\mathrm{Hom}}

\newcommand{\Loc}{\mathrm{Loc}}

\numberwithin{equation}{subsection}

\usepackage{tikz}

\newcommand{\SF}[1]{S^1_{\bF_{#1}}}

\newcommand{\GL}{\mathrm{GL}}
 
\title{The Alexander polynomial at prime powers}
\author{David Treumann}

\begin{document}
\maketitle

\begin{abstract}
When $K$ is a knot and $p \gg 0$ is a prime, we discuss a finite set whose cardinality is $\Delta_K(p^n)$, the value of the Alexander polynomial of $K$ at $p^n$.
\end{abstract}

\section{Introduction}

This paper and its sequel explores a literal interpretation of the standard analogy between the Frobenius in the absolute Galois group of a finite field, and the generator of the fundamental group of a circle.

We define an ``$F$-field'' on a manifold $M$ to be a locally constant sheaf of algebraically closed fields of positive characteristic.  The best example is the $F$-field $\underline{k}$ on a circle whose fiber at a base point is an algebraic closure $k$ of $\bF_p$, and whose monodromy around the circle is the $p$th power map --- we write $\SF{p}$ to indicate the circle along with such an $F$-field.  Every other good example is pulled back from this one along a map $\mathfrak{f}:M \to \SF{p}$.

The pullback of $\underline{k}$ along $\mathfrak{f}$ is a sheaf of rings, and we can consider sheaves of modules over it.  For example, the category of locally free $\underline{k}$-modules of finite rank on $\SF{p}$ is equivalent to the category of finite-dimensional $\bF_p$-vector spaces.  In this paper we will prove the following:
\begin{prop*}
If $M$ is the complement of a knot in $S^3$, $p$ is a sufficiently large prime, and $\mathfrak{f}:M \to \SF{p}$ has degree $n$ on $H_1(M)$, then the set of isomorphism classes of invertible $\mathfrak{f}^*\underline{k}$-modules is a finite commutative group of cardinality $\Delta_K(p^n)$.
\end{prop*}

In the Proposition, $\Delta_K(x) = c_0 + c_1 x + \cdots + c_n x^n$ denotes the classical Alexander polynomial of $K$, with Alexander's normalization that $c_0$ should be a positive integer.  The primes $p$ that must be excluded are the divisors of $c_0$.

The Proposition gives a combinatorial interpretation of the Alexander polynomial that seems to be new, though it can be proved by unwinding the sheaf-theoretic definitions until one arrives at the cokernel of Alexander's matrix, evaluated at $x = p^n$.  We will express this ``unwinding of definitions'' as a brief development of theory in \S\ref{sec:two}, and give the proof in \S\ref{sec:three}.

A sequel paper will make a more detailed study of the influence of $F$-fields on categories of constructible sheaves and other categories associated to symplectic manifolds.  The Proposition perhaps gives some motivation for doing so.

In contrast with knot complements and most other 3-manifolds, when $M$ is two-dimensional the locally constant sheaves of $\mathfrak{f}^*\underline{k}$-modules come in positive-dimensional moduli over the ground field $k$.  Local systems over a nontrivial $F$-field are much more stable than traditional local systems, in the GIT sense that their automorphism groups are finite.  When $M$ is symplectic, one can study a Fukaya category whose objects are Lagrangian submanifolds $L \subset M$ equipped with locally constant sheaves of $(\mathfrak{f}^*\underline{k})\vert_L$-modules.  The sequel will develop some examples related to Deligne-Lusztig varieties.

\subsection{Notation}
\label{intro:notation}
Throughout, $p$ is a prime number.  We let $\bF_p$ denote the finite field with $p$ elements, and fix an algebraic closure $k$ of $\bF_p$ once and for all.  If $q$ is a power of $p$, we write $\bF_q$ for the subfield of $k$ with $q$ elements, fixed by the $q$th-power automorphism of $k$.

\section{$F$-fields and local systems}
\label{sec:two}

For each prime power $q$ we fix an oriented circle $\SF{q}$, endowed with a base point.  For definiteness, we put $\SF{q} = \bR/\log(q) \bZ$, oriented in the positive direction and with the base point at the coset of $0$.  We call these the ``reference circles.''

We endow each reference circle with a sheaf of rings $\underline{k} = \underline{k}_{\SF{q}}$.  It is locally constant with fiber $k$ (as in \S\ref{intro:notation}), and the monodromy from $0$ to $\log(q)$ is given by $a \mapsto a^q$.  More formally, define the \'etal\'e space of $\underline{k}$ to be the quotient of $\bR \times k$ by the equivalence relation
\begin{equation}
\label{eq:espace-etale}
\underline{k} = (\bR \times k) / \sim \qquad (t,x) \sim (t + \log(q),x^q)
\end{equation}

\subsection{$F$-fields}
\label{subsec:2.1}
We define an \emph{$F$-field} of characteristic $p$ on a space $X$ to be a continuous map $\mathfrak{f}_X:X \to \SF{p}$.  When $q = p^\nu$ is a power of $p$ we have a canonical $F$-field $\SF{q} \to \SF{p}$ sending $t+ \log(q) \bZ$ to $t + \log(p)\bZ$.

When $X$ carries an $F$-field $\mathfrak{f}:X \to \SF{p}$, it carries a locally constant sheaf of rings $\mathfrak{f}^* \underline{k}$ as well, the pullback of \eqref{eq:espace-etale} along $\mathfrak{f}$.  If the $F$-field is clear from context we will sometimes write $\underline{k}_X$ or just $\underline{k}$ in place of $\mathfrak{f}^* \underline{k}$. 
We write $\Loc(X,\mathfrak{f}^* \underline{k})$ for the category of locally free sheaves of $\mathfrak{f}^* \underline{k}$-modules on $X$.  It is an abelian category but it is not usually $k$-linear. For example, $\Loc(\SF{q},\underline{k})$ is equivalent to the category of finite-dimensional $\bF_q$-modules.  The equivalence is given by the global sections functor $L \mapsto \Gamma(\SF{q},L)$, whose $\bF_q$-module structure comes from $\Gamma(\SF{q},\underline{k}) \cong \bF_q$.

\subsection{Variants of $\underline{k}$}
Any automorphism at all can be repurposed as a locally constant sheaf on a circle; e.g. we could replace $k$ with any perfect ring $R$ along with its $p$th power automorphism.  If we let $\underline{R}$ denote the corresponding sheaf of rings on $\SF{p}$, then studying $\Loc_n(X,\mathfrak{f}^* \underline{R})$ as $R$ runs through perfect $k$-algebras gives a moduli functor that is represented by a perfect Deligne-Mumford stack.  We will study these and similar moduli problems in the sequel paper --- in the case of a knot complement these moduli stacks are zero-dimensional.  Here are three other significant examples of a different flavor:

\subsubsection{Isocrystals} 
Let $\bQ_{p^{\infty}}$ be the maximal unramified extension of $\bQ_p$, whose residue field is $k$.  The Frobenius automorphism of $k$ lifts to a field automorphism of $\bQ_{p^{\infty}}$ that we denote by $\sigma$.  Then $\bQ_{p^{\infty}}$ is the fiber of a sheaf of rings on $\SF{p}$ whose monodromy is $\sigma$, we denote it by $\underline{\bQ_{p^{\infty}}}$.  A locally constant sheaf of $\underline{\bQ_{p^{\infty}}}$-modules is the same data as an isocrystal, studied in \cite{Dieudonne}.  

\subsubsection{The Tate motive}
\label{subsec:Tate}
The action of $\pi_1(\SF{p})$ on the multiplicative group $k^*$ of $k$ corresponds to a local system of abelian groups that we denote by $\underline{k^*}$.  There is a closely related local system of $\bZ[1/p]$-modules, which we denote by $\bZ[1/p](1)$ --- the fiber is $\bZ[1/p]$ and the generator of $\pi_1(\SF{p})$ acts by multiplication by $p$.  The two sheaves are related by 
\[
\text{(sheaf hom)}(\bZ[1/p](1),\text{const. sheaf with fiber }k^*) \cong \underline{k^*}
\]

\subsubsection{Finite Chevalley groups}
\label{subsec:Chevalley}
If $\bG$ is an algebraic group over $k$, and $\sigma:\bG \to \bG$ is the Frobenius isogeny coming from an $\bF_q$-rational structure on $\bG$, then $\sigma$ induces an automorphism on $k$-points $\bG(k) \to \bG(k)$.  We denote the corresponding locally constant sheaf of groups on $\SF{q}$ by $\underline{(\bG(k),\sigma)}$.  
The groupoid of $\underline{(\bG(k),\sigma)}$-torsors over $\SF{q}$ is equivalent to the classifying groupoid of torsors over the finite group $\bG(k)^{\sigma}$.

There is a similar construction for the Suzuki or Ree isogenies of $\bG = \mathrm{Sp}_4$ when $p = 2$ and $\bG = \mathrm{G}_2$ when $p = 3$ --- the square of one of these isogenies is a Frobenius map for an $\bF_q$-rational structure.  Thus they again induce bijections on $k$-points and it is natural to regard $\underline{(\bG,\sigma)}$ as a sheaf of groups on a reference circle of circumference $\log(\sqrt{q})$.

\subsection{Framings} 
We write $\Loc_n(X,\mathfrak{f}^* \underline{k})$ for the groupoid of locally free sheaves of rank $n$ in $\Loc(X,\mathfrak{f}^*\underline{k})^{\simeq}$ and all isomorphisms between them.  There are two ways to present the groupoid as a quotient:
\subsubsection{Point framings} Let $x \in X$ be a base point and let $E \in \Loc_n(X,\mathfrak{f}^*\underline{k})$.  A \emph{point framing} of $E$ at $x$ is a $(\mathfrak{f}^* \underline{k})_x$-basis in $E_x$.  Note $(\mathfrak{f}^* \underline{k})_x$ is canonically identified with $\underline{k}_{\mathfrak{f}(x)}$.  The set of isomorphism classes of point-framed local systems of rank $n$ is in bijection with the set of homomorphisms
\begin{equation}
\label{eq:point-framing}
\rho:\pi_1(X,x) \to \GL_n(\underline{k}_{\mathfrak{f}(x)}) \rtimes \pi_1(\SF{p})  
\end{equation}
that commute with the projections to $\pi_1(\SF{p})$.  In the semidirect product in \eqref{eq:point-framing}, the distinguished generator of $\pi_1(\SF{p})$ acts on $\GL_n(\underline{k}_{\mathfrak{f}(x)})$ by raising each matrix entry to the $p$th power.

If we write $\rho(\gamma) = \rho_1(\gamma) \rtimes \mathfrak{f}(\gamma)$, the bijection sends $\rho$ to the locally constant sheaf of abelian groups whose fiber above $x$ is $\underline{k}_{X,x}^n$, and whose monodromy around the loop $\gamma$ acts on the vector $(v_1,\ldots,v_n)$ by
\[
\rho_1(\gamma)(v_1^{p^{\mathfrak{f}(\gamma)}},\ldots,v_n^{p^{\mathfrak{f}(\gamma)}})
\]
The groupoid $\Loc_n(X,\mathfrak{f}^*\underline{k})$ is equivalent to the quotient of the set of \eqref{eq:point-framing} by the $\GL_n(\underline{k}_{\mathfrak{f}(x)})$-conjugation action.

\subsubsection{Section framings}
With $q = p^\nu$, we shall call the data of a map $s:\SF{q} \to X$ commuting with the projections to $\SF{p}$ a \emph{$\nu$-sheeted multisection} of the $F$-field $\mathfrak{f}$.  Each $E \in \Loc_n(X,\mathfrak{f}^*\underline{k})$ restricts to a locally constant sheaf on such a multisection, or equivalently an $n$-dimensional $\bF_q$-vector space, $\Gamma(s,E)$.  A \emph{section framing} of $E$ at $x$ is an $\bF_q$-basis for this vector space.

A section framing induces a point framing, at the image of the base point of $\SF{q}$ under $s$.  Using $g$ to denote the distinguished generator of $\pi_1(\SF{p})$, we say that  homomorphism \eqref{eq:point-framing} preserves the section-framing if it carries $s$ to $1 \rtimes g^{\nu}$.  These homomorphisms are in bijection with the set of section-framed local systems.  A change of basis in $\Gamma(s,E)$ corresponds to conjugating $\rho$ by an element in the finite group $\GL_n(\underline{k}_{\bF_{p^\nu}}) \subset \GL_n(\underline{k}_{\mathfrak{f}(x)})$.  In particular this shows that $\Loc_n(X,\mathfrak{f}^*\underline{k})$ has finite isotropy groups whenever $X$ is connected and the $F$-field is nontrivial (i.e. whenever $\pi_1(X) \to \pi_1(\SF{p})$ is nonzero).

\subsection{Invertible modules}

When $n = 1$, the homomorphisms \eqref{eq:point-framing} do not necessarily factor through an abelian quotient of $\pi_1(X,x)$.  Nevertheless the tensor product over $\mathfrak{f}^*\underline{k}$ endows the set of isomorphism classes in $\Loc_1(X,\mathfrak{f}^*\underline{k})$ with the structure of commutative group. The groupoid $\Loc_1(X,\mathfrak{f}^*\underline{k})$ has a symmetric monoidal structure --- it is a commutative $2$-group.  

Indeed, if we regard the abelian group $\GL_1(\underline{k}_{\mathfrak{f}(x)})$ as a $\pi_1(X,x)$-module through the homomorphism $\pi_1(X,x) \to \pi_1(\SF{p})$, and $\rho_1:\pi_1(X,x) \to \GL_1(\underline{k}_{\mathfrak{f}(x)})$ is a $1$-cocycle on $\pi_1(X,x)$ with coefficients in this module, then $\rho_1 \rtimes \mathfrak{f}$ is a homomorphism of the form \eqref{eq:point-framing}.  This is a bijection between such cocycles and point-framed local systems. If we write $Z^1 := Z^1(\pi_1(X,x),\GL_1(\underline{k}_{\mathfrak{f}(x)}))$ for this group of cocycles, the commutative $2$-group structure on $\Loc_1(X,\mathfrak{f}^*\underline{k})$ can be encoded by the two-term chain complex
\begin{equation}
\label{eq:GL1Z1}
\GL_1(\underline{k}_{\mathfrak{f}(x)}) \to Z^1
\end{equation}
where the differential is the usual differential in group cohomology.  If $\pi_1(X,x) \to \pi_1(\SF{p})$ is nonzero, the kernel is a finite group (it is $\GL_1$ of a finite subfield of $\underline{k}_{\mathfrak{f}(x)}$).

\section{Proof of the Proposition}
\label{sec:three}

Fix a diagram $D$ for a knot with $v$ crossing points and $v+2$ regions.  Let us label the regions $0,\ldots,v+1$, and orient $K$.  Then there is a Dehn presentation of $\pi_1(S^3 - K)$ with $v+1$ generators and $v$ relations \cite[\S 8]{Alexander}.  The generators $g_i$ correspond to regions $1,\ldots,v+1$ of $D$, and we put $g_0 = 1$.  Each crossing gives a relation between the generators associated to the four regions incident with it, as in the following diagram,
\begin{center}
\begin{tikzpicture}
\node at (-3,0) {$g_j g_k^{-1} g_{\ell} g_m^{-1} = 1$};
\draw[thick] (-1,0)--(1,0);
\draw[thick] (0,-1)--(0,-.2);
\draw[thick, ->] (0,.2)--(0,1);
\node at (-.25,-.25) {$\cdot$};
\node at (-.25,.25) {$\cdot$}; 
\node at (.5,.5) {$m$};
\node at (-.5,.5) {$j$};
\node at (-.5,-.5) {$k$};
\node at (.5,-.5) {$\ell$};
\end{tikzpicture}
\end{center}
(In \cite{Alexander}, knot diagrams are drawn with two dots on the left side of the underpass crossing --- we have put them in the diagram above as well.)

The ``index'' of a region defined by Alexander \cite[Fig. 2]{Alexander} gives a homomorphism 
\[
I:\pi_1(S^3 - K) \to \bZ
\]  
A choice of prime $p$ and prime power $q = p^\nu$ turns the index homomorphism into an $F$-field $\mathfrak{f}$, by identifying it with a homotopy class of maps $S^3 - K \to \SF{q}$.  A rank one point-framed local system of $\mathfrak{f}^* \underline{k}$-modules on $S^3 -K$ is completely specified by a family of scalars $z_j \in \GL_1(k)$, one for each region. In terms of these scalars, the action of $g_j$ (and its inverse, recorded for convenience) on $x \in k$ is given by
\[
g_j(x) = x^{q^{I(j)}} z_j \qquad \text{(and }g_j^{-1}(x) = x^{q^{-I(j)}} z_j^{-q^{-I(j)}} \text{)}
\]
The scalars $z_j$ must obey $z_0 = 1$ and an additional relation for each crossing indicident with regions $j,k,\ell,m$ as above, which reduces to 
\[
z_j z_k^{-q} z_{\ell}^{q} z_m^{-1} = 1 \qquad \text{ or } \qquad z_j z_k^{-q^{-1}} z_{\ell}^{q^{-1}} z_m^{-1} = 1
\]
according to whether the crossing is left-handed or right-handed, respectively.  
An element $y \in \GL_1(k)$ acts on the point-framing by sending $(z_j)_{j = 1}^{v+1}$ to $(y^{I(j) - 1} z_j)_{j = 1}^{v+1}$.  Thus, the group of isomorphism classes of objects in $\Loc_1(X,\mathfrak{f}^* \underline{k})$ is isomorphic to the middle cohomology of the chain complex
\begin{equation}
\label{eq:vacuum}
\GL_1(k) \to \GL_1(k)^{v+1} \to \GL_1(k)^{v}
\end{equation}
where the first and second differentials are
$
y \mapsto (y^{q^{I(j)} - 1})_{j = 1}^{v+1}$ and  $(z_j)_{j = 1}^{v+1} \mapsto (z_j z_k^{-q^{\pm 1}} z_{\ell}^{q^{\pm 1}} z_m^{-1})_c$.  

Consider the $v \times (v+1)$-matrix whose rows are indexed by the crossings, whose columns are indexed by non-null regions, and whose $(c,j)$-entry is indicated by the following diagram if $j$ is incident with $c$, and is otherwise $0$:
\begin{center}
\begin{tikzpicture}
\draw[thick, ->] (-1,0)--(1,0);
\draw[thick] (0,-1)--(0,-.2);
\draw[thick, ->] (0,.2)--(0,1);
\node at (-.25,-.25) {$\cdot$};
\node at (-.25,.25) {$\cdot$}; 
\node at (.5,.6) {$-1$};
\node at (-.5,.6) {$1$};
\node at (-.7,-.6) {$-x$};
\node at (.65,-.6) {$x$};
\end{tikzpicture}
\qquad
\begin{tikzpicture}
\draw[thick,<-] (-1,0)--(1,0);
\draw[thick] (0,-1)--(0,-.2);
\draw[thick, ->] (0,.2)--(0,1);
\node at (-.25,-.25) {$\cdot$};
\node at (-.25,.25) {$\cdot$}; 
\node at (.5,.6) {$-1$};
\node at (-.5,.6) {$1$};
\node at (-.75,-.6) {$-x^{-1}$};
\node at (.75,-.6) {$x^{-1}$};
\end{tikzpicture}
\end{center}
Let $D'$ denote the diagram obtained from $D$ by switching the sense of over and under at every crossing.  If we multiply each row corrsponding to a right-hand crossing by $x$, we obtain the usual Alexander matrix (i.e. the coefficient matrix of the system of equations \cite[Eq. 3.3]{Alexander}) for the diagram $D'$, with the column corresponding to the null region left off.  In particular $A(x)$ is elementary equivalent to the usual Alexander matrix for both $D'$ and (using the mirror invariance of the Alexander invariants) $D$.

By evaluating $A(x)$ at $x = q = p^\nu$, we obtain a matrix $A(q) \in \bZ[1/p]^{v \times (v+1)}$.  Then \eqref{eq:vacuum} is obtained by taking $\Hom(C,\GL_1(k))$, where $C$ is a complex of free $\bZ[1/p]$-modules of the form
\[
\bZ[1/p]^v \xrightarrow{A(q)} \bZ[1/p]^{v+1} \xrightarrow{} \bZ[1/p]
\]

In particular, the order of the middle cohomology of \eqref{eq:vacuum} is equal to the order of the middle cohomology of this complex, which is equal to the order of the torsion subgroup of the cokernel of $A(q)$.  This in turn is equal to the prime-to-$p$ part of the greatest common divisor of all the $v \times v$-minors of $A(q)$.  If $p$ does not divide the constant term of $\Delta_K(x)$, that divisor is just $\Delta_K(q)$.

\section{Plausible generalizations}
It is natural to consider the size of $\Loc_1(M;\mathfrak{f}^*\underline{k})$ ``as an orbifold'', i.e. to weight each isomorphism class by the reciprocal of the order of its automorphism group.  In the case of a knot complement these automorphism groups are not very sensitive to $\mathfrak{f}$ so that the orbifold count is $\Delta_K(q)/(q-1)$.  In the case of a more general $3$-manifold both the kernel and cokernel of \eqref{eq:GL1Z1} are more irregular as $\mathfrak{f}$ varies, but the orbifold count is likely to have a clean relation to the multivariable Alexander polynomial.

One can study local systems of modules over $\mathfrak{f}^* \underline{\bQ_{p^{\infty}}}$.  By restricting such a local system $L$ to a meridian of the knot, one gets a rank one isocrystal whose slope is a discrete invariant of $L$.  The set of $L$ of a fixed slope make a $p$-adic analytic manifold (in fact a torsor for a commutative $p$-adic analytic group) whose dimension over $\bQ_{p^n}$ is the degree of the Alexander polynomial.  I suspect that these $p$-adic manifolds carry natural measures, perhaps up to powers of $p$, of volume $\Delta_K(p^n)$.

One might obtain interesting invariants by counting nonabelian local systems of $\mathfrak{f}^* \underline{k}$-modules, or more generally torsors for the sheaves of  groups $\mathfrak{f}^* (\underline{\bG},\sigma)$.  The problem of doing so can be expressed as the problem of solving high-degree equations in the entries of a matrix in $\bG(k)$.  For example (using the Wirtinger presentation of $\pi_1(S^3 - K)$ in place of the Dehn presentation) a  $\mathfrak{f}^*\underline{(\bG,\sigma)}$-torsor on the complement of a trefoil is determined by an element $g \in \bG(k)$ subject to the equation
\[
g \sigma^2(g) = \sigma(g)
\] 
and taken up to the conjugation action by the finite group $\bG^{\sigma}$.  But even in the simplest examples I have not been able to solve these equations directly when $\bG$ is not commutative.

\subsection*{Acknowledgments}
I am grateful to the Institute for Advanced Study where this paper was written.  My stay at the IAS was supported by a Sloan fellowship, a von Neumann fellowship, and a Boston College faculty fellowship, and I was also supported by NSF-DMS-1510444.

\end{document}